\documentclass{article}

\begin{document}

\title{Algebraic families of constant scalar curvature K\"ahler metrics}
\author{Simon Donaldson}
\newcommand{\uX}{\underline{X}}
\newcommand{\cX}{{\cal X}}
\newcommand{\cL}{{\cal L}}
\newcommand{\cY}{{\cal Y}}
\newcommand{\cZ}{{\cal Z}}
\newcommand{\Chow}{{\rm Chow}}
\newcommand{\Ch}{{\rm Ch}}
\newcommand{\cO}{{\cal O}}
\newcommand{\oV}{\overline{V}}
\newcommand{\tFut}{\widetilde{{\rm Fut}}}
\newcommand{\Fut}{{\rm Fut}}
\newcommand{\bC}{{\bf C}}
\newcommand{\bP}{{\bf P}}
\newcommand{\dbd}{i\partial \overline{\partial}}
\newcommand{\cdbd}{\partial\overline{\partial}}

\newcommand{\db}{\overline{\partial}}
\newcommand{\Vol}{{\rm Vol}}
\newcommand{\Riem}{{\rm Riem}}
\newcommand{\bZ}{{\bf Z}}
\newcommand{\tomega}{\tilde{\omega}}
\newcommand{\tDelta}{\tilde{\Delta}}
\newcommand{\uP}{\underline{P}}
\newcommand{\Aut}{{\rm Aut}}
\newcommand{\Tr}{{\rm Tr}}
\newcommand{\Osc}{{\rm Osc}}
\newcommand{\cU}{{\cal U}}
\newcommand{\cV}{{\cal V}}
\newcommand{\oJ}{\overline{J}}
\newtheorem{defn}{Definition}
\newtheorem{prop}{Proposition}
\newtheorem{lem}{Lemma}
\newtheorem{cor}{Corollary}
\newtheorem{thm}{Theorem}

\maketitle



\section{Introduction}

Let $\uX\rightarrow S$ be a family of Fano manifolds of complex dimension $n$, where the base $S$ is a quasi-projective variety. Thus for each $s\in S$ we have an $n$-dimensional Fano manifold $X_{s}$. Define $S^{*}\subset S$ to be the set of parameters $s$ such that $X_{s}$ admits a K\"ahler-Einstein metric. Then  we will prove:
\begin{thm} Suppose that for each $s\in S$ the automorphism group of $X_{s}$ is finite. Then $S^{*}$ is a Zariski-open subset of $S$.
\end{thm}
Of course, the statement allows the possibility that $S^{*}$ is empty. It is well-known, and easy to prove, that under these hypotheses $S^{*}$ is open in the $C^{\infty}$-topology; the theorem can be viewed as stating that the condition of admitting a
 K\"ahler-Einstein metric  is of an algebro-geometric nature, without specifying what this is. 
The assumption on the automorphism group cannot be removed, as shown by the
example of the Mukai-Umemura manifold \cite{kn:T}, \cite{kn:D2}. It seems
likely that the same result holds if the identity component of the automorphism
group of each $X_{s}$ is isomorphic to some fixed Lie group (or, equivalently,
if the dimension of the space of holomorphic vector fields is constant),
but  some technical obstacles arise in extending the arguments to this case.

We should say straightaway that this is not a new result, it was proved by Odaka in \cite{kn:Odaka2}. Odaka's proof relies on the solution of Yau's conjecture, relating K\"ahler-Einstein metrics and stability by Chen, Donaldson and Sun \cite{kn:CDS} which does, at least in principle, specify the criterion for existence in algebro-geometric terms. One reason for giving this second proof is that it is independent of \cite{kn:CDS}. The Zariski-open statement is easier to prove than the result in \cite{kn:CDS} so it seems worth having a direct proof. The verification of the stability condition is not really practical in the current state of the art, so the simpler statement about Zariski-openness seems at present almost as useful as the more detailed one.  Further, while there is overlap, the
approach here is substantially different to that in \cite{kn:CDS} and we hope that the techniques have interest and may be relevant to other problems. This approach was discussed in an unpublished document  \cite{kn:D3}, dating from 2009 and the previous paper \cite{kn:D4}. There was a central difficulty, discussed in those references. which it turns out can be overcome using the subsequent results of Sun and the author in \cite{kn:DS}.  See also the further discussion in 2.2 below.

One difference between the techniques used here and in \cite{kn:CDS} is that the latter is more tied to the (second order) K\"ahler-Einstein equation while our techniques focus on the (fourth order) constant scalar curvature condition. In this direction, we will prove a more general result. Suppose now that the quasi-projective variety $S$ parametrises a family $X_{s}$ of polarised manifolds. So there is a line bundle on the total space $\uX$ which restricts to a positive line bundle $L_{s}$ on each fibre $X_{s}$. Now we let $S^{*}$ be  defined by the condition that $X_{s}$ admits a constant scalar curvature K\"ahler metric $\omega_{s}$ in the class $c_{1}(L_{s})$. In such a case the isometry class of metric $\omega_{s}$ is unique.
\begin{defn}
    The family $\uX\rightarrow S$ has Property R if there is a constant $C$ such that for all $s\in S^{*}$ the diameter and the $L^{\infty}$ norm of the Ricci curvature of $(X_{s},\omega_{s})$ are bounded by $C$.
\end{defn}

Then we have
\begin{thm} Suppose that $S$ has property R and that for each $s\in S$ the automorphism group of $X_{s}$
is finite. Then $S^{*}$ is a Zariski-open subset of $S$.\end{thm}

Theorem 2 implies Theorem 1 (the Ricci curvature bound being trivial, and the diameter bound following from Myers' Theorem), and it is Theorem 2 which we will prove in this article. As things stand, the generalisation has only theoretical interest, since we have no way of establishing Property R beyond the K\"ahler-Einstein case. But the more general statement emphasises the point that many of our arguments focus on the constant scalar curvature
 condition.  
 
 In Section 2 below we will explain the proof of Theorem 2, assuming various more technical statements which are established in the later sections of the paper.
There are two  main new features in  our approach. On the differential-geometric side, our basic method is to combine arguments involving the Bergman function and the \lq\lq volume estimate'' established by the X. Chen and the author \cite{kn:CD}, and by Cheeger and Naber \cite{kn:CheegerNaber}. On the algebro-geometric side, we work with general degenerations, not assuming a $\bC^{*}$-action. This is  in may ways a more natural setting, and offers some advantages which can be useful in other contexts.

The author is grateful to Xiuxiong Chen, Yuji Odaka and Song Sun for discussions related to this paper. The author was partially supported by ERC Advanced Investigator Award 247331.

\section{Proof of Theorem 2}
\subsection{Degenerations}
We begin with a simple observation.
\begin{lem}
  To prove Theorem 2 it suffices to prove that, for any family $S$ satisfying the hypotheses of that Theorem, if $S^{*}$ is non-empty then there is  a non-empty Zariski-open subset of $S$ which is contained in $S^{*}$.
\end{lem}

This is straightforward. In the context of Theorem 2, if $S^{*}$ is empty then the statement is trivially true. Suppose $S^{*}$ is not empty and thus, following the statement of the lemma, contains a non-empty Zariski-open subset $S\setminus Z$. Let $S_{1} $ be a component of maximal dimension of $Z$. Then we can take $S_{1}$ in place of $S$ in the discussion and we have a subset $S_{1}^{*}= S^{*}\cap S_{1}\subset S_{1}$ corresponding to constant scalar curvature metrics. Either $S_{1}^{*}$ is empty, in which case $S\setminus S^{*}$ contains $S_{1}$, or non-empty. In the latter case, the complement
$S_{1}\setminus S_{1}^{*}$ is, by hypothesis, contained in a lower-dimensional algebraic set. Proceeding in this way over all the components, and by induction on dimension, one sees that $S^{*}$ is Zariski-open.

We now come to a central topic of this paper.
\begin{defn}
A degeneration of $(X,L)$ is a flat analytic family $\pi:\cX\rightarrow \Delta$ over the disc $\Delta\subset \bC$ and a line bundle $\cL\rightarrow \cX$ which is very ample   on every fibre of $\pi$, together with an isomorphism between the restriction of $\cL\rightarrow \cX$ to the 
 punctured disc $\Delta^{*}$ and the product family $(X,L^{m})\times \Delta^{*}$, for some integer $m\geq 1$.
\end{defn}

We call the integer $m$ the {\it multiplicity} of the degeneration. Soon we will see that we can reduce to the case when $m=1$, but for the moment we keep this parameter. We will usually just refer to a degeneration by its total space $\cX$ etc. A particular case of a degeneration is an {\it equivariant degeneration}, or \lq\lq test configuration'', when the family extends to $\bC$ and there is a compatible $\bC^{*}$ action on $\cX$. We need a notion of {\it non-trivial} degenerations. If we work with general schemes this is a slightly issue, as pointed out by Li and Xu \cite{kn:LX}, but we will only need to consider degenerations with normal central fibres and for these non-triviality just means that the degeneration is not equivalent to a product. 

We want to define a notion of a manifold $X_{s}$, for $s\in S$, having \lq\lq generic degenerations at multiplicity $m$''. To begin, we can clearly suppose that the line bundle $L_{s}\rightarrow X_{s}$ is very ample for {\it some} $s\in S$ and replacing $S$ by a non-empty Zariski open subset we can then suppose this is true for {\it all} $s\in S$. Thus for each $m\geq 1$ and each $s\in S$ the sections of $L_{s}^{m}\rightarrow X_{s}$ give a projective embedding of $X_{s}$. Similarly (by passing to a Zariski open subset if necessary) we can suppose that all higher cohomology groups vanish so that the dimension of $H^{0}(X_{s}, L_{s}^{m})$ is a fixed number $N(m)+1$ for all $s\in S$, given by the Riemann-Roch formula. If we choose a basis for the space of sections,  the image of the embedding corresponds to  point in the Chow variety $\Chow=\Chow(N(m), d(m))$ of cycles of degree $d(m)$ in $\bP^{N(m)}$, where $d(m)= m^{n} c_{1}(L_{s})^{n}$. Let $\Gamma_{s}$ be the $SL(N_{m}+1)$-orbit of this point in $\Chow$; in other words, the set of points we get from any choice of basis. Up to a possible covering of the disc, a degeneration of $X_{s}$ of multiplicity $m$ corresponds to a holomorphic map $f:\Delta\rightarrow
\Chow$ with $f(t)\in \Gamma_{s}$ for $t\in \Delta^{*}$. If we have a degeneration $\cX$ we get a map $f$ when we choose a
 trivialisation of the bundle $\pi_{*}(\cO(L^{m}))$. Conversely, if we have such a map $f$ we pull back the universal family to the disc and after possibly taking a covering (in the case when the automorphism group is non-trivial) we can trivialise the family over the punctured disc.

We also want to consider a notion of a family of degenerations (of multiplicity $m$). Let $N$ be an open set in $S$. We restrict the family $\uX$ to $N$ and then lift to the product $N\times \Delta^{*}$ to get a space
 $\cU_{N}\rightarrow N\times \Delta^{*}$. By a family of degenerations parametrised by $N$ we mean a flat family $\cV\rightarrow N\times \Delta$ with an isomorphism
$\Phi$ between the restriction of $\cV$ to $N\times \Delta^{*}$ and $\cU_{N}$. We also require a line bundle over $\cV$ with the obvious conditions, so the upshot is that for each $s$ in $N$ the restriction of $\cV$ to $\{s\} \times \Delta$ gives a degeneration of $X_{s}$ of multiplicity $m$. 

Define $ J\subset S\times \Chow $ to be the pairs $(s,Z)$ with $Z\in \Gamma_{s}$ and let $\oJ$ be its closure in $S\times \Chow$. 
\begin{defn}
We define a subset $GD_{m}\subset S$. A point $\sigma$ is in $GD_{m}$ if for each point $(\sigma, [Y])]$ in $\oJ$  the following two conditions hold.
\begin{itemize}
\item $[Y]$ lies in the closure $\overline{\Gamma_{\sigma}}$ of the orbit corresponding to $X_{\sigma}$;
\item There is neighbourhood $N$ of $\sigma$ in $S$ and a family of degenerations
$\cV\rightarrow N\times \Delta$ of multiplicity $m$, parametrised by $N$, with fibre $Y$over $(\sigma,0)$.
\end{itemize}\end{defn}

With this definition in place we can state one of our more technical results, which is proved in Section 5.
\begin{prop}
For each $m$ there is a non-empty, Zariski-open subset $\Omega_{m}$ contained in $GD_{m}\subset S$.
\end{prop}

Applying Lemma 1 we get the following corollary:
\begin{cor}
To prove Theorem 2 it suffices to show that there is an $m$ such that, if $s_{i}$ is a sequence in $S^{*}$ converging to a point $\sigma$ in $GD_{m}$  then $\sigma$ also lies in $S^{*}$.\end{cor}

All of the above is  preparatory to the main proof 
which is now to establish the statement in Corollary 1. 

We digress here slightly to explain the significance of the subset $GD_{m}$ of \lq\lq generic degenerations''. The essential condition is the first one, which rules out the phenomenon of \lq\lq splitting of orbits''. Consider, as a simple example, the action of $\bC^{*}$ on $\bC^{2}$ with weights $(1,-1)$. For $s\in \bC$ let $\Gamma_{s}$ be the $\bC^{*}$-orbit of $(1,s)$ and let $J$ be the set of pairs $(s,Z)$ where $s\in \bC$ and $Z\in \Gamma_{s}\subset \bC^{2}$. Let  $\sigma=0$ and $Z_{0}= (0,1)\in \bC^{2}$. Then $(\sigma, Z_{0})$ lies in the closure of $J$ but  $Z_{0}$ does not lie in the closure of the orbit $\Gamma_{\sigma}$. In the context of Corollary 1 one would like to establish a different statement, to wit that if $s_{i}\in S^{*}$ converge to $\sigma$ and if $X_{\sigma}$ is {\it K-stable} then $\sigma\in S^{*}$. This would imply that if $S^{*}$ is not empty then it consists precisely of the  K-stable points. The splitting of orbits is a fundamental difficulty here, since one has to deal with the possibility that the constant scalar curvature metrics $\omega_{s_{i}}$ might converge to a metric on another variety $X'$ but where $X'$ is not the central fibre of any degeneration of $X_{\sigma}$. We refer to the further discussion in \cite{kn:D3}. For our purposes here the point is that if we are satisfied with the \lq\lq Zariski open'' statement we can avoid this problem.

\subsection{Limits}

We take up the situation considered in Corollary 1, with a sequence $s_{i}$ in $S^{*}$ converging to $\sigma \in S$. We will write $X_{i}$ for $X_{s_{i}}$ and $L$ for $L_{s_{i}}$.  Fix $m\geq 1$ for the moment. For each $i$ we have an $L^{2}$ norm on sections of $L^{m}$ over $X_{i}$ and we define embeddings $T_{i,m}:X_{i}\rightarrow \bC\bP^{N(m)}$ using orthonormal bases with respect to these norms. So the $T_{i,m}$ are uniquely defined up to the action of the compact group
$U(N(m)+1)$. Taking a subsequence, we can suppose that the images $T_{i,m}(X_{i})$ converge to some limiting algebraic cycle $W_{m}\subset \bP^{N(m)}$. Then we have two alternatives:
\begin{enumerate}\item $W_{m}$ is the image of an embedding of $X_{\sigma}$ by sections of $L^{m}$. In terms of the Chow variety, the corresponding point
$[W_{m}]$ lies in the orbit $\Gamma_{\sigma}$. 
\item $W_{m}$ is not the image of an embedding of $X_{\infty}$
by sections of $L^{m}$. In terms of the Chow variety, the corresponding point
$(\sigma, [W_{m}])$ lies in $\oJ$ but not in the orbit $\Gamma_{\sigma}$.
\end{enumerate}

\begin{prop}
We can choose an $m_{0}$ so that, for all such sequences $s_{i}$ in $S$, if $m\geq m_{0}$ and
alternative
(1) above holds then $X_{\sigma}$ admits a constant scalar curvature metric, i.e. $\sigma$ lies in $S^{*}$.
\end{prop}
The proof, which is relatively standard,  is given in Section 4 below.

Given this result, we can focus attention on alternative (2). Thus we have algebraic cycles $W_{m}$ for $m\geq m_{0}$, not equivalent to $X_{\sigma}$. Here we reach a crucial point in the discussion. There is {\it a priori} the possibility that the $W_{m}$ are essentially different for different values of $m$. This could only happen if they fail to be normal. This situation does not fit in well with the theory of degenerations and their numerical invariants discussed in the following subsection. The  notion of 
\lq\lq b-stability'' introduced in \cite{kn:D4} was designed   to address precisely this difficulty. However it turns out this is not necessary in our problem, since we have the following result.
\begin{prop}
We can choose $m_{1}\geq m_{0}$ so that $W_{m_{1}}$ is a normal projective variety and for all $k\geq 1$ the natural map from $W_{km_{1}}$ to $W_{m_{1}}$ is an isomorphism.
\end{prop}
This is a part of what is proved in \cite{kn:DS}. Replacing $L$ by $L^{m_{1}}$
we may as well suppose that $m_{1}=1$ and we just write $W$ for the limiting variety. We also write $GD\subset S$ for $GD_{1}$.
 To sum up our arguments so far; to prove Theorem 2 it suffices to show that  if the point 
$\sigma$ lies in $GD\subset S$ then the second alternative above cannot hold

{\bf Remark} The notion of b-stability arises from  generalisations of the Chow and Futaki invariants discussed below to allow for a varying family $W_{m}$.  It seems quite possible that this approach could be pushed through to prove Theorem 2 without use of Proposition 3 but, at best,  this would need a good deal of technical work.

\subsection{Chow and Futaki invariants}

We return to consider a degeneration $\cX$ of $(X,L)$ and from now on we can take multiplicity $m=1$. Recall that we have a line bundle $\cL\rightarrow \cX$ so for each $k\geq 1$ we can take the direct image of the sheaf of holomorphic sections of $\cL^{k}$ which we just write as $\pi_{*}(\cL^{k})$. This is a holomorphic vector bundle over the disc $\Delta$. The trivialisation of the family $\cX$ over the punctured disc $\Delta^{*}$ induces a trivialisation of this vector bundle there. Relative to this trivialisation, we have a first Chern class, which can be viewed as an integer. We define:
\begin{equation} \tau_{k}(\cX)= -c_{1}(\pi_{*}(\cL^{k}))\in \bZ. \end{equation}

\begin{lem}
For large enough  $k$ the integer $\tau_{k}(\cX)$ is a polynomial function of $k$ of degree at most $(n+1)$ \end{lem}
To see this we can use the trivialisation of the family $\cX$ over $\Delta^{*}$ to extend to a family over $\bP^{1}$. Then  apply the
 Grothendieck-Riemann-Roch theorem.

Define $I=I(\cX)$ to be the co-efficient of $k^{n+1}$ in the polynomial $k\mapsto \tau_{k}(\cX)$. Let $N(k)+1$ be the dimension of $H^{0}(X, L^{k})$, as before, and write $$V= \frac{1}{n!} c_{1}(L)^{n}$$
so that $N(k)\sim V k^{n}$ as $k\rightarrow \infty$, by Riemann-Roch. For each $k$, we define the Chow invariant to be
\begin{equation} \Ch_{k}(\cX)= \frac{\tau_{k}(\cX)}{N(k)}- k \frac{I(\cX)}{V}. \end{equation}
By construction, the numbers $\Ch_{k}(\cX)$ have  a finite limit as $k$ tends to infinity and we define the {\it Futaki invariant}:
\begin{equation} \tFut(\cX)= \lim_{k\rightarrow \infty} \Ch_{k}(\cX) \end{equation}

The heart of the proof of Theorem 2 is contained in two statements about this Futaki invariant.
\begin{prop}
Let $\cX$ be a degeneration with central fibre $W$, where $W=\lim T_{i,k}(X_{i})$ is the limiting variety of subsection (2.2). Then $\tFut(\cX)\leq 0$.
\end{prop}
This represents the differential-geometric input to Theorem 2, using the fact that  the $X_{i}$ have metrics of constant scalar curvature. The proof of Proposition 4 is given in Section 3.

\begin{prop}
Let $\cX$ be a non-trivial degeneration of $(X,L)$ where $X$ has a metric of constant scalar curvature in the class $c_{1}(L)$ and ${\rm Aut}(X,L)$ is finite. Then $\tFut(\cX)>0$.
\end{prop}

This is a result of Szekelyhidi \cite{kn:Gabor}, extending that of Stoppa \cite{kn:Stoppa}. We discuss  this, and give another proof, in Section 5.

With all this in place, we can now give the proof of Theorem 2. We want to show that if $\sigma$ is in $GD\subset S$ then $[W]$
is in the orbit $\Gamma_{\sigma}$. We suppose not and derive a contradiction.
Since $(\sigma, [W])$ lies in $\oJ$ we can apply the definition of $GD$
and we have a  family of degenerations $\cU\rightarrow N\times \Delta$ with fibre $W$ over $(\sigma,0)$. We can suppose that $N$ is connected. For $s\in N$ we restrict $\cU$ to $\{s\}\times \Delta$ to get a degeneration $\cX_{s}$ of $X_{s}$. Proposition 4 shows that $\tFut(\cX_{\sigma})\leq 0$. For large enough $i$ we have $s_{i}\in N$ and
Proposition 5 shows that $\tFut(\cX_{s_{i}})>0$. But it is clear from the definition that $\tFut(\cX_{s})$ is a constant function of $s\in N$ so this gives the desired contradiction.

\section{Proof of Proposition 4}
\subsection{The moment map}
We review some general theory, as discussed also in \cite{kn:D4}.  Let $\Chow$ be the Chow variety whose points represent $n$-dimensional algebraic cycles of degree $n! V$ in $\bP^{N}$. The Chow construction gives an embedding of $\Chow$ in a projective space $\bP(U)$ and the action of $GL(N+1)$ on $\Chow$ is induced from  a representation of $GL(N+1)$ on $U$. Write $\lambda_{\Chow}$ for the hyperplane bundle restricted to $\Chow$. This can also be described as a Deligne pairing, see the exposition in \cite{kn:PS}. Let $\cX$ be a degeneration which can be embedded in $\bP^{N}\times \Delta$ with $n$-dimensional fibres of degree $d$ and with line bundle $\cL\rightarrow \cX$ the lift of $\cO(1)$ to the product. As in Section 2.1, we get a map $f:\Delta\rightarrow \Chow$ and the pull-back of $\lambda_{\Chow}$ is a line bundle over $\Delta$. We can write $f$ as 
$$   f(t)= g(t) [X]$$ where $g$ is a holomorphic map from $\Delta^{*}$ to $GL(N+1)$ which extends to a meromorphic matrix-valued function over the disc. This map $g$ gives a trivialisation of $f^{*}(\lambda_{\Chow})$ over the punctured disc, so we have in integer $c_{1}(\lambda_{\Chow})$, the first Chern class of the line bundle relative to this trivialisation. The basic fact, going back to Knudsen and Mumford \cite{kn:KM}, is that this is another description of the Chow invariant in that we have
\begin{equation} I(\cX)= - c_{1}(\lambda_{\Chow}). \end{equation}
Explicitly, if  $\widehat{[X]}$ is a point in $U$ lying over $[X]$ then when $g$ takes values in $SL(N+1)$ the Chow invariant is $V^{-1}$ times the order of the pole of the meromorphic vector-valued function $g(t)\widehat{[X]}$ at $0$. The condition that $g$ takes values in $SL(N+1)$ corresponds to the vanishing of $\tau=\tau_{1}$, in the notation of 2.3 above.  More generally, if $g$ maps  $\Delta^{*}$ to
 $GL(n+1,\bC)$, let $a$ be the order of the pole of $g(t) [X]$ and $b$ be the order of the pole of $\det g$. Then the Chow invariant is
\begin{equation} \Ch(\cX)= \frac{a}{V}- \frac{b}{N+1}. \end{equation}

Now fix the standard Hermitian metric on $\bC^{N+1}$. This gives a Fubini-Study metric on $\bP^{N}$ and hence a volume form $d\mu_{FS}$ on any algebraic cycle. For an $n$-dimensional cycle $Z$ of degree $n! V$  we define a matrix
\begin{equation}   M(Z)_{\alpha \beta}= V^{-1} \left[ \int_{Z} \frac{z_{\alpha}\overline{z}_{\beta}}{\vert z\vert^{2}} d\mu_{FS} \right]_{0} \end{equation}
where $[\ ]_{0}$ denotes projection to the trace-free part. 
This gives a map
$$   \sqrt{-1}M: \Chow\rightarrow {\rm Lie}(SU(N+1)), $$
and the basic fact is that this is a {\it moment map} for the isometric action of $SU(n+1)$
on $\Chow$ with respect to a certain metric. (We pass over the point that
$\Chow$ may be singular, since inspection of the arguments shows that this is not in the end relevant.) If $\cX$ is an equivariant degeneration with respect to a $\bC^{*}$-action generated by $A\in \sqrt{-1} {\rm Lie}(SU(N+1))$ then we have a standard formula,
$$  \Ch(\cX)= \langle M(W), A \rangle, $$
where $W$ is the central fibre. The crucial point for us here is that we can extend this to an {\it inequality} for any degeneration $\cX$. Recall that a map $g(t)$ as above can be factorised as:
  \begin{equation}  g(t) = L(t) t^{A} R(t) \end{equation}
  where $L, R$ are holomorphic functions on the disc taking values in $SL(N+1)$ and $A$ is Hermitian with integer eigenvalues (so that $t^{A}$ is well-defined).  Then we have
 
\begin{equation}   \Ch(\cX)\leq \langle M(W), A \rangle \end{equation}
where $W$ is the central fibre. This is proved in \cite{kn:D4}, Lemma 1.  

We return to the situation considered in 2.2 above, so we have a sequence $X_{i}$ of manifolds with constant scalar curvature metrics $\omega_{i}$ and for each $k$ we have embeddings $T_{k,i}:X_{i}\rightarrow \bP^{N(k)}$ whose images converge to $T_{k}(W)$. Clearly
$$M(T_{k}(W))=\lim_{i\rightarrow \infty} M(T_{k,i}(X_{i}). $$
Let  $\cX$ be any degeneration with central fibre $W$. Replacing $L$ by $L^{k}$, we have a Hermitian endomorphism $A_{k}$ as above and (8) gives
\begin{equation}  \Ch_{k}(\cX) \leq \lim_{i\rightarrow \infty} \langle M( T_{k,i}(X_{i})),A_{k}\rangle. \end{equation}
For Hermitian matrices $A,M$ write $N(A)$ for the difference between the largest and smallest eigenvalues of $A$  and write $\vert M\vert_{1}$ for the sum of the absolute values of the eigenvalues of $M$ (the \lq\lq trace norm''). Then if $M$ has trace zero we have
  \begin{equation}   \vert \langle A, M\rangle\vert \leq 2 N(A)\vert M\vert_{1} \end{equation}
 To prove Proposition 4 we need two facts. The first is purely algebraic:
\begin{lem}  For any degeneration $\cX$  \begin{equation}N(A_{k})= k N(A_{1}) \end{equation} \end{lem}
The second contains all the differential geometry:
\begin{prop}
 Assuming that $S$ has property $R$, there is a constant $C_{S}$ such that
 \begin{equation} \vert M(T_{k}(X_{s})) \vert_{1} \leq C_{S}  k^{-2} \log k \end{equation}
 for all $s\in S^{*}$. \end{prop}
 For, assuming these statements, we get from (9),(10), (11), (12) that
 $$\Ch_{k}(\cX)\leq 2 C_{S} C_{\cX} k^{-1} \log k, $$
 so, passing to the limit, we have $\tFut(\cX)\leq 0$.

\subsection{Application of the volume estimate}
In this subsection we  write $X, L$ for some $X_{s}, L_{s}$ with $s\in S_{*}$. We write
 $\omega$ for the constant scalar curvature metric in the class $c_{1}(L)$ on  $X$ and we fix a Hermitian
metric on $L$ with corresponding curvature form $-2\pi i\omega$. We write
$\omega_{FS,k}$ for the pull-back by $T_{k}$ of the Fubini-Study metric and
$$d\mu_{FS,k}= k^{-n} (n!)^{-1}\omega_{FS,k}^{n}. $$
Notice that we have normalised so that the integral of the volume $d\mu_{FS,k}$ over $X$ is the fixed number $V$, independent of $k$. 

Recall the \lq\lq density of states'' function $\rho_{k}$ on $X$ which can be defined as follows. Let $s_{\alpha}$ be any orthonormal basis of the holomorphic sections of $L^{k}$ with respect to the $L^{2}$ norm. (This $L^{2}$ norm uses the induced fibrewise metric on $L^{k}$ and the volume form $k^{n} (n!)^{-1} \omega^{n}$.) Then
\begin{equation}  \rho_{k} = \sum_{\alpha} \vert s_{\alpha}\vert^{2}. \end{equation}
Let, as above, $A$ be any Hermitian endomorphism of $H^{0}(X,L^{k})$ and choose $s_{\alpha}$ to be a orthonormal basis of eigenfunctions, with eigenvalues
$\lambda_{\alpha}$. Now define a function $H_{A}$ 
by
$$   H_{A}= \rho_{k}^{-1} \left( \sum_{\alpha} \lambda_{\alpha} \vert s_{\alpha}\vert^{2}\right).
$$
 Then, unravelling the definitions, we find that
  \begin{equation} \langle M_{k}(X),A\rangle = V^{-1} \int_{X} H_{A}\
   d\mu_{FS,k}. \end{equation}
   
   By duality, the estimate stated in Proposition 4 is equivalent to showing that for all $A$ and $k$  we have
   \begin{equation}   \vert \langle M_{k}(X),  A \rangle\vert \leq C_{S} k^{-2}\log k \vert A \vert_{\infty},\end{equation}
where $\vert A\vert_{\infty}=\max_{\alpha} \vert \lambda_{\alpha}\vert$. Since $M_{k}(X)$ has trace zero it suffices to prove this when $A$ has trace zero. 
Clearly the $L^{\infty}$ norm of $H_{A}$ is bounded by $\vert A\vert_{\infty}$ and if $A$ has zero trace then 
$$  \int_{X} H_{A} \rho_{k} d\mu_{FS,k}=0. $$
So to prove Proposition 4  it suffices  to prove that for {\it any} function $H$ with $$  \int_{X} H \rho_{k} d\mu_{FS,k} =0$$ we have
\begin{equation} \vert \int_{X} H d\mu_{FS,k}\vert\leq V C_{S} \Vert H\Vert_{L^{\infty}} k^{-2} \log k. \end{equation}

In our situation we may assume that  the dimension of $H^{0}(X,L^{k}$ is given by 
Riemann-Roch, so:
$$  {\rm dim} H^{0}(X,L^{k})= V k^{n} P(k) $$ where
$$  P(k)= 1+ \alpha_{1} k^{-1} + \dots \alpha_{n} k^{-n}, $$
for fixed numbers $\alpha_{i}$. To prove Proposition 4 it suffices to show that we have a
bound
\begin{equation} \vert \int_{X} H \left( \frac{\rho}{P(k)} \omega^{n}- d\mu_{FS,k}
\right) \vert \leq V C_{S} k^{-2} \log k \Vert H\Vert_{L^{\infty}}\end{equation}
for {\it all} functions $H$. Finally, if we define a $2n$-form $\Theta_{k}$ by
$$  \Theta_{k}= \frac{\rho_{k}}{P(k)} \omega^{n} - \omega_{FS,k}^{n} , $$
it suffices to prove that 
\begin{equation} \int_{X} \vert \Theta_{k} \vert \leq V^{-1} C_{S} k^{-2}\log k. \end{equation}
From the definition of the Fubini-Study metric we find that
$$  \omega_{FS,k}= \omega+ k^{-1} \dbd \log \rho_{k} , $$
so we can express everything in terms of $\rho_{k},P_{k}$:
\begin{equation}  \Theta_{k} = \frac{\rho_{k}}{P(k)} \omega^{n} - (\omega+
k^{-1} \dbd \log \rho_{k})^{n}. \end{equation}

The crucial ingredient in the proof is the asymptotic expansion of $\rho_{k}$
as $k\rightarrow \infty$. By well-known results of Tian, Caitlin, Lu, Zelditch,
Ruan and others, for any {\it fixed} metric $\omega$ we have an asymptotic
expansion
$$  \rho_{k} \sim (2\pi)^{-n}( 1+ a_{1} k^{-1} + a_{2} k^{-2} + \dots) , $$
where the $a_{i}$ are functions on $X$ determined by the curvature tensor
and, in particular, $a_{1}$ is one half of the scalar curvature. If
 the metric has constant scalar curvature, as we are assuming, $a_{1}$
is a constant and the fact that the integral of $\rho_{k}$ is $k^{n}$
times the dimension of $H^{0}(X,L^{k})$ implies that this constant coincides
with the leading co-efficient $\alpha_{1}$ of $P(k)$.

For our purposes we need an extension of the standard results which can be
applied uniformly to {\it any} of our metrics. Given $r>0$ we define $Z_{r}\subset
X$ to be the $r$-neighbourhood of the points where $\vert \Riem \vert\geq
r^{-2}$ (all definitions being with respect to the metric $\omega$).
Let $\Omega_{r}= X\setminus Z_{r}$ be the complement. Thus for a point $x\in
\Omega_{r}$ the curvature is bounded by $r^{-2}$ on the ball of radius $r$
about $x$. 
\begin{prop}
Write $\rho_{k}= (2\pi)^{-n}(1+ \alpha_{1} k^{-1} + \eta_{k})$.
Under the assumption that $S$ has property R there are constants 
$b_{0},C_{0}, C_{1}, C_{2}$ such that at a point $x\in \Omega_{r}$ and for $k\geq b_{0}^{2} r^{-2}$ we have
$$  \vert \eta_{k}\vert \leq C_{0} k^{-2} r^{-4},\ \vert \nabla \eta_{k}\vert
\leq C_{1} k^{-2} r^{-5}, \ \vert \nabla\nabla \eta_{k}\vert \leq C_{2} k^{-2}
r^{-6}. $$
\end{prop}
We discuss the proof in Section 6 below.
 
Next we have a bound on the volume of $Z_{r}$.
\begin{prop}
Under the assumption that $S$ has Property R there is a constant $C_{3}$
such that ${\rm Vol} (Z_{r}) \leq C_{3} r^{4}$.
\end{prop}

This is proved in \cite{kn:CD}. More precisely, the statement there is for K\"ahler-Einstein metrics but it is straightforward to extend the proof to constant scalar curvature K\"ahler metrics,  in a fixed integral cohomology class,  with  bounds on the Ricci curvature and diameter. See  the remarks in Section 4 below.

Now we go on to complete the proof of Proposition 6.

Choose $b\geq \max( (4C_{0})^{1/4}, b_{0})$ and write $r_{0}= b k^{-1/2}$. Set $\Omega=\Omega_{r_{0}},
Z=X\setminus \Omega$. Thus at points of $\Omega$ we have $\vert \eta_{k}\vert
\leq 1/4$. We can assume that $k$ is large enough that $\vert\alpha_{1} k^{-1}\vert\leq
1/4$ and so $\vert \rho-1\vert \leq 1/2$.  We can then estimate $\dbd \log
\rho_{k}$ in an elementary way, from the statements in Proposition 7,  to
get
  $$  \vert \dbd \log \rho \vert \leq C_{5}  k^{-2} r^{-6}$$ in $\Omega_{r}$,
for any $r\geq r_{0}$. Then
$$  \vert k^{-1} \dbd \log \rho \vert \leq C_{5} b^{2} k^{-2} r^{-4}.$$
We can suppose $b$ chosen large enough that $\vert k^{-1} \dbd \log \rho\vert
\leq 1/2$ in $\Omega$. Write $d\mu$ for the volume form of the metric $\omega$ and $\omega_{FS}^{n}= (1+F_{1}) d\mu$.  By Proposition 7 we have
$$  \vert F_{1} \vert \leq C_{6} k^{-2} r^{-4}$$ in $\Omega_{r}$, for $r\geq
r_{0}$ . Likewise if we write $\rho/P= 1+F_{2}$ then, from the first item
in Proposition 7, we have $$ \vert F_{2}\vert \leq C_{7} k^{-2} r^{-4}$$
in $\Omega_{r}$, for $r\geq r_{0}$. 

 By definition $$ \int_{X} \vert \Theta \vert  = \int_{X} (1+F_{1})-(1+F_{2}) d\mu $$
so $$  \int_{X}\vert \Theta \vert \leq \int_{\Omega} \vert F_{1} \vert +
\vert F_{2} \vert\  d\mu + \int_{Z} \vert 1+F_{1} \vert + \vert 1+F_{2} \vert. \ d\mu . $$
The definition of $F_{1}$ implies that $1+F_{1}>0$. We can suppose $k$ is
large enough that $P_{k}>0$ thus $\rho_{k}/P_{k}>0$ and hence also $1+F_{2}>0$.
Thus
$$\int_{X}\vert \Theta \vert \leq \int_{\Omega} \vert F_{1}\vert +\vert F_{2}\vert
d\mu + \int_{Z} 2+ F_{1}+F_{2}\ d\mu.  $$
Since the integral of $\omega_{FS}^{n}$ over $X$ is equal to that of $\omega^{n}$
we have $\int_{X} F_{1} d\mu = 0$. Similarly the fact the integral of $\rho$ is
$P_{k}$ times the volume of $X$ implies that $\int_{X} F_{2}d\mu = 0$. Thus
$$  \int_{Z} F_{1}+F_{2}\  d\mu = -\int_{\Omega} F_{1}+F_{2} \ d\mu \leq \int_{\Omega} \vert
F_{1} \vert + \vert F_{2} \vert\  d\mu . $$
So we get
\begin{equation}  \int_{X} \vert \Theta\vert \leq 2\int_{\Omega} \vert F_{1}
\vert +\vert F_{2} \vert \ d\mu + 2 \Vol(Z). \end{equation}
Our volume estimate of Proposition 8  gives $\Vol(Z)\leq C_{4} r_{0}^{4}= C_{4} b^{4} k^{-2}$.
To handle the first term in(20), we have $\vert F_{1}\vert +\vert F_{2}\vert\leq
C_{8} k^{-2} r^{-4}$ say  in $\Omega_{r}$, for $r\geq r_{0}$. Write
$ \vert F_{1} \vert + \vert F_{2} \vert = C_{8} k^{-2} F$ so that $F\leq
r^{-4}$ in $\Omega_{r}$, for $r\geq r_{0}$. Let $f(t)$ be the distribution
function of $F$ on $\Omega$
$$  f(t)= \Vol\{ x\in \Omega: F(x)\geq t\}. $$
Thus, for $r\geq r_{0}$, $f(r^{-4})$ is bounded by the volume of $Z_{r}$
and hence by $C_{4} r^{4}$. In other words $f(t)\leq C_{4} t^{-1} $ for $t\leq
r_{0}^{-4}$ and $f(t)$ vanishes for $t\geq r_{0}^{-4}$. Clearly also we have
$f(t)\leq \Vol(X)$ for all $t$. Then
$$  \int_{\Omega} F =\int_{0}^{r_{0}^{-4}} f(t) dt\leq \int_{0}^{r_{0}^{-4}}
\min (\Vol(X), C_{4} t^{-1})\  dt. $$
Let $\tau= C_{4}/\Vol(X)$. We may suppose that $k$ is large enough that
$r_{0}^{-4}> \tau$. Then we get
$$\int_{\Omega} F \leq C_{4}+ C_{4} \int_{\tau}^{r_{0}^{-4}} t^{-1} dt, $$
which is bounded by $C_{9}+ C_{4} \log r_{0}^{-4}$. Then
$$\int_{\Omega} \vert F_{1}\vert +\vert F_{2} \vert \leq C_{8} k^{-2} (C_{9}+
C_{4} \log ( b^{-4} k^{2}). $$
Combining with the estimate for the last term in (20) we get the required bound
$$  \int_{X}\vert \Theta \vert \leq C_{10} k^{-2} \log k. $$

\section{Regularity theory: proof of Proposition 2}

Recall that we have a sequence $s_{i}\rightarrow \sigma$ in $S$ and constant scalar curvature metrics $\omega_{i}$ on $X_{i}=X_{s_{i}}$. Choose a family of reference metrics $\tomega_{i}$ on $X_{i}$, converging in $C^{\infty}$ to a limit $\tomega_{\sigma}$ (with respect to a local $C^{\infty}$ trivialisation of the family).  For each $i$ we can write
$$ \omega_{i}- \tomega_{i}= \sqrt{-1}\cdbd \phi_{i}. $$

The first step is to prove that, under our hypotheses (bounded diameter and Ricci curvature) the oscillation
$$ {\rm Osc}(\phi_{i})= {\rm max}\ \phi_{i}- {\rm min}\ \phi_{i} $$
satisfies a fixed bound. This can be read off from the results of \cite{kn:DS} but we prefer to minimise our dependence on \cite{kn:DS} here so we recall another well-known argument (which has the advantage of being effective, given an effective bound in the volume estimate of Proposition 8).

  Let $T_{i}, \tilde{T_{i}}: X_{i}\rightarrow  \bP^{N}$ be the embeddings defined by the metrics
$\omega_{i}, \tomega_{i}$. By Proposition 7 (and replacing $L$ by a fixed power $L^{k}$ if necessary), we can suppose that there is an $\omega_{i}$-ball $B_{i}\subset X_{i}$ of a fixed radius on which the metric $\omega_{i}$ is well-approximated by the pull back $T_{i}^{*}(\omega_{FS})$ of the Fubini-Study metric under $T_{i}$. The hypothesis of Proposition 2 is that the limit of $T_{i}(X_{i})$ is in the same orbit and this 
means that $T_{i}= g_{i} \circ \tilde{T}_{i}$ where $g_{i}$ is a bounded sequence of projective transformations. Thus the pull-backs $T_{i}^{*}(\omega_{FS}), \tilde{T}_{i}^{*}(\omega_{FS})$ are quasi-isometric and it follows that on $B_{i}$ the metrics $\tilde{\omega_{i}}, \omega_{i}$ are quasi-isometric. Note that the diameter and Ricci bounds for the metrics $\omega_{i}$ give a lower bound on the volume of $B_{i}$. 

Now we have a standard result.
\begin{lem}
Suppose $\omega,\tomega$ are two K\"ahler metrics on $X$, each with diameter and $\vert Ricci\vert$ bounded by $C$. Suppose there is an $\eta>0$ and an open set $B\subset X$ such that $\tomega\leq \eta^{-1} \omega$ on $B$ and the $\omega$-volume of  $B$ exceeds $\eta$. Then if $\omega=\tomega+\dbd \phi$ there is a bound $\Osc(\phi)\leq K$, where $K$ depends only on $n,C,\eta$.
\end{lem} 

We recall the proof. Set $f=\phi-\max \phi$. Since $\tomega+i\dbd f >0$ we have $\tDelta f \geq -n$ where $\tDelta$ is the Laplacian of the metric $\tomega$ (with the \lq\lq analysts'' sign convention). If $x_{0}$ is a point where $\phi$ attains its maximum we have
$$   0 = f(x_{0})= {\rm Av}_{\tomega} f - \int_{X} G(x_{0}, y) \tDelta f_{y} d\tilde{\mu}_{y}, $$
(in obvious notation). We normalise the Green's function so that $G(x_{0},y)\geq 0$ for all $y$. It is then known that 
$$  \int_{X} G(x_{0}, y) d\tilde{\mu}_{y} \leq \kappa, $$
where $\kappa$ depends only on $C,n$. (This fact can be extracted from \cite{kn:CL} for example.)
So we get $$\int_{X} f  d\tilde{\mu} \geq - n \kappa \Vol(X).$$
In terms of $\phi$:
$$  \int_{X} (\max \phi - \phi)\  d\tilde{\mu} \leq n \kappa \Vol(X), $$ so
\begin{equation}  \int_{B} (\max \phi -\phi)\  d\tilde{\mu} \leq n \kappa \Vol(X). \end{equation}
Interchanging the roles of  $\omega, \tilde{\omega}$
$$  \int_{B} (\phi -\min \phi)\ \ d\mu \leq n \kappa \Vol(X). $$
Thus \begin{equation} \int_{B} (\phi-\min \phi)\ d\tilde{\mu} \leq \eta^{-n} n \kappa \Vol(X), \end{equation}
since $d\tilde{\mu}\leq \eta^{-n} d\mu$.  Adding (21) and (22),
$$  \int_{B} (\max\phi-\min \phi)\  d\tilde{\mu}= (\max \phi -\min \phi) \Vol(B,d\tilde{\mu})
 \leq (1+ \eta^{-n}) n \kappa \Vol(X) $$
So $$ \max\phi -\min \phi \leq \eta^{-1}(1+\eta^{-n}) n \kappa \Vol(X). $$

In our case we apply this using the balls $B_{i}$ to get a bound on $\Osc \phi_{i}$. Once we have this bound on the K\"ahler potential, the assumed bound on the Ricci curvature allows us to appeal to Theorem 5.1 in \cite{kn:CH} which
gives   upper and lower bounds on the metrics
\begin{equation} K_{1}^{-1} \tomega_{i}\leq \omega_{i}
\leq K_{1} \tomega_{i}. \end{equation}
(This can also be achieved by combining the Yau estimate with the Chen-Lu estimate: there is an   good exposition of the latter in \cite{kn:JMR}.) To get higher estimates we can proceed as follows. Let 
$$L_{i}= \log \frac{\omega_{i}^{n}}{\tomega_{i}^{n}}. $$
Then $\cdbd L_{i}= \sqrt{-1}(\rho_{i}- \tilde{\rho}_{i})$ where $\rho_{i}, \tilde{\rho_{i}}$ are the Ricci $(1,1)$-forms. The Ricci curvature bound and (23) imply that $\tDelta_{i} L_{i}$ are bounded in $L^{\infty}$ and hence $L_{i}$ are bounded in $C^{1,\alpha}$. Then the estimates of Cafarelli \cite{kn:Caf}, for the Monge-Amp\`ere operator, using again (23), give a $C^{1,\alpha}$ bound on the metrics $\omega_{i}$. From there on we can use the constant scalar curvature equation to get estimates on all higher derivatives.

\

To conclude this section we take up the discussion from (3.2), of the extension of the volume estimate from \cite{kn:CD} to our situation (constant scalar curvature, bounded Ricci curvature, bounded diameter). Examining the proofs in \cite{kn:CD}, one sees that what is needed is the following.
\begin{lem} Suppose that $\Omega_{i}$ is a sequence of constant scalar curvature K\"ahler metrics with bounded Ricci curvature and satisfying a uniform volume non-collapsing condition. Suppose there is a pointed Gromov-Hausdorff limit $Z,\Omega_{\infty}$ and that $R\subset Z$ is the regular set. Then the $\Omega_{i}$ converge in $C^{\infty}$ over $R$.
\end{lem}
  
  The standard Anderson-Cheeger-Colding theory for non-collapsed metrics with bounded Ricci curvature gives local
   $C^{1,\alpha}$ convergence of the metrics in harmonic co-ordinates, so $C^{,\alpha}$ convergence of the Christoffel symbols.  The K\"ahler condition states that the complex structures $J_{i}$ are covariant constant, so we get convergence to a $C^{1,\alpha}$ limiting almost-complex structure $J_{\infty}$. The Newlander-Nirenberg theorem extends to the $C^{1,\alpha}$ case (\cite{kn:NW}) and this gives  $C^{1,\alpha}$ convergence of the metrics
in local holomorphic coordinates. Then the constant scalar curvature condition gives $C^{\infty}$ convergence, as above.

\section{Algebraic geometry}

\subsection{Review of theory: results of Stoppa and Szekelyhidi}

The main goal of this section is to establish Proposition 5. This extends a result of Stoppa \cite{kn:Stoppa}, in the case of $\bC^{*}$-equivariant degenerations. The complete argument, in that case, involves a number of steps which we will now review.

Recall the notion of \lq\lq bi-asymptotic stability'' defined in \cite{kn:D4}.
For a point $x$ in $X$, write $\hat{X}_{x}$ for the blow-up and $L_{E}\rightarrow
\hat{X}_{x}$ for the line bundle defined by the exceptional divisor. For
sufficiently large  $\mu$, and then for $r$ sufficiently large, the line bundle $L^{r\mu}\otimes
L_{E}^{r}$ is  ample on $\hat{X}_{x}$. We say that $X$ is $(r,\mu)$-stable
if, for all $x$, the image of $\hat{X}_{x}$ under this projective embedding
is Chow stable. We say that $(X,L)$ is bi-asymptotically stable if there
is a $\mu_{0}$ and for $\mu\geq \mu_{0}$ an $r_{0}(\mu)$ such
that $X$ is $(r,\mu)$ stable for $r\geq r_{0}(\mu)$. 

Now we have the following fact.
\begin{prop}
    Suppose that $X$ admits a constant scalar curvature K\"ahler metric in the class $c_{1}(L)$ and $\Aut(X,L)$ is finite. Then $(X,L)$ is bi-asymptotically stable.
\end{prop}
The proof of this combines two ingredients
\begin{enumerate}
\item According to Arezzo and Pacard \cite{kn:AP}, there is an $\epsilon(X,L)>0$ such that any blow-up $\hat{X}_{x}$ admits a constant scalar curvature K\"ahler metric in the class $c_{1}(L)+ \epsilon c_{1}(L_{E})$  for $0<\epsilon<\epsilon(X,L)$.
\item Let $(Y,\Lambda)$ be any polarised manifold. According to the result of \cite{kn:D1}, for large enough $k$ the projective embedding of $Y$ by sections of $L^{k}$ can be chosen to satisfy the \lq\lq balanced'', or zero-moment map, condition $M(Y)=0$. From the interpretation of this in terms of Chow theory going back to Zhang \cite{kn:Zhang}, it follows that the image is Chow stable.
\end{enumerate}
Then the Proposition follows, taking $Y= \hat{X}_{x}$ and 
$\mu_{0}>\epsilon(X,L)^{-1}$. Both of the steps above use the hypothesis that $\Aut(X,L)$ is finite. 

The central result of Stoppa's paper \cite{kn:Stoppa} is an algebro-geometric statement---if $(X,L)$ is bi-asymptotically stable then for any non-trivial $\bC^{*}$-equivariant degeneration $\cX$ with generic fibre $(X,L)$ we have $\Fut(\cX)>0$. To prove this he considers  the blow-up $\widehat{\cX}$ of $\cX$ along a  certain subscheme $\Sigma\subset \cX$, which yields a degeneration of the blow-up of $X$ at a point. The bi-asymptotic stability  hypothesis implies that,  for a suitable line bundle over $\widehat{\cX}$, we have $\Fut(\widehat{\cX})\geq 0$ and, by careful calculation, Stoppa shows that $\Fut(\widehat{\cX})<\Fut(\cX)$.

Turning back to our situation, to prove Proposition 6 it suffices to establish the algebro-geometric result:
\begin{prop}
    If $(X,L)$ is bi-asymptotically stable and $\Aut(X,L)$ is finite then for any non-trivial degeneration $\cX$ of $(X,L)$ we have $\tFut(\cX)>0$.
\end{prop}

In \cite{kn:Gabor}, Szekelyhidi developed a wide-ranging extension of the notion of stability of a pair $(X,L)$ involving filtrations of the ring 
$\bigoplus_{k} H^{0}(X,L^{k})$ and proved an analogue of Stoppa's result. He also showed how to associate a filtration to  a degeneration $\cX$ and then Proposition 10 follows as a special case of his main result.  However it seems worthwhile to have a proof tailored to this particular situation, avoiding some of the technicalities in Szekelyhidi's  work, and this is what we will take up in 5.3 below.

\subsection{Some elementary results about degenerations}

We consider a degeneration $\cX\subset \bP(V)\times \Delta$, so the fibres $X_{t}$ can be viewed as varieties in $\bP(V)$ and the line bundle $\cL$ is the pull-back of $\cO(1)$. In our situation we can assume that the fibres do not lie in any proper projective subspace of $\bP(V)$. As we have discussed in (3.1), this can be represented by a meromorphic matrix-valued function $g:\Delta^{*} \rightarrow GL(V)$ and this function can be factorised as $g(t)= L(t) t^{A} R(t)$. By applying an automorphism of $\bP(V)\times \Delta$---that is to say, changing the trivialisation of $\pi_{*}\cL$---we can suppose that $L(t)=1$ (with the same $R(t))$. We can also suppose that $R(0)=1$.  Choose a basis so that  that $A$ is diagonal with eigenvalues $\lambda_{0}\geq \lambda_{1}  \dots\geq \lambda_{\min}$.
\begin{lem} We can choose $g$ so that $R(t)$ is lower triangular and 
the entry $R_{\alpha \beta}$ is a polynomial of degree less than
 $\lambda_{\alpha}-\lambda_{\beta}$.
\end{lem}

 To see this we can first factorise $R(t)= Q(t) \tilde{R}(t)$ where $Q$ is upper triangular and $\tilde{R}$ is lower triangular. Then $L(t)= t^{A} Q(t) t^{-A}$ is holomorphic across $t=0$ and 
$$  t^{A} R(t)= L(t) t^{A} \tilde{R}(t). $$
Applying an automorphism, as above, we can replace $R(t)$ by $\tilde{R}(t)$. So we can suppose that in fact $R(t)$ is lower triangular. Now it is easy to see that we can factorise a lower triangular matrix function $R(t)$ as
$$   R(t)= R_{1}(t) R_{2}(t), $$
where $R_{1}, R_{2}$ are lower triangular and each $(\alpha \beta)$ entry of $R_{2}$ is a polynomial of degree less than $\lambda_{\alpha}-\lambda_{\beta}$ and each $(\alpha \beta)$ entry of $R_{1}$ contains terms of order at least
$\lambda_{\alpha}-\lambda_{\beta}$. This condition on $R_{1}$ means that
$L(t)= t^{A} R_{1}(t) t^{-A}$ is holomorphic across $0$ and
$$  t^{A} R(t)= L(t) t^{A} R_{1}(t). $$
 
 As before we can now suppose that $R=R_{1}$, which has the form stated in the Lemma.

 It is convenient to normalise $g$ so that the largest eigenvalue $\lambda_{0}$
is $0$. This is just achieved by changing the trivialisation of the line bundle over $\Delta^{*}$ by multiplying by a power of $t$, which does not affect the Chow or Futaki invariants.   With this choice we will
 say that $\cX$ is a {\it normalised degeneration}. 
 \begin{lem}
 For a non-trivial normalised degeneration $\cX$ we have $I(\cX)<0$.
 \end{lem}
 
 The compactified family $\overline{\cX}$ is embedded in a projective bundle
$\bP(E)$ where
$E$ is a vector bundle over $\bP^{1}$ defined by the transition function
$g$. Explicitly we have a trivialisation $V\times \Delta_{\infty}$
over a neighbourhood of infinity and a trivialisation $V \times \Delta$
over a neighbourhood of $0$, and a point $(x,t)$ in the first trivialisation
is identified with a point $(y,t)$ in the second where $y=g_{t}(x)$. Let
$\cL\rightarrow \bP(E)$ be the dual of the tautological line bundle. Then
-$\tau_{k}$ is the degree of the vector bundle 
$\pi_{*}(\cL\vert_{\overline{\cX}})$ over $\bP^{1}$. By the 
Grothendieck-Riemann-Roch
theorem, the leading term  in the expansion of -$\tau_{k}$ is given by $1/n!$ times the evaluation of evaluating
$c_{1}(\cL)^{n+1}$ on $\overline{\cX}$. This leading term is $-I$ by definition, so  the Lemma will  follow if we
establish that the sections of $\cL$   separate
generic points of $\overline{\cX}$. To see this, we consider the co-ordinates
$y_{\alpha}$ as sections of $\pi_{*}(\cL)$ over a neighbourhood of $0$ and $x_{\alpha}$
as sections over a neighbourhood of infinity. They are related by
   $$  y_{\alpha}= t^{\lambda_{\alpha}} x_{\alpha} + \sum_{\beta<\alpha} t^{-\lambda_{\alpha}} R_{\alpha \beta}
x_{b} . $$
Since the terms in $R_{ab}$ have degree less than $\lambda_{\alpha}-\lambda_{\beta}$
we see that $y_{\alpha}$ extends holomorphically over infinity.  Further, if $\lambda_{\alpha}>0$ then $t y_{\alpha}$ also extends holomorphically over infinity. Since the fibres
of $\overline{\cX}$ do not line in any proper projective subspace it is clear
that the sections separate generic points of $\overline{\cX}$.

We now move on to the proof of Lemma 3. This will follow from an intrinsic characterisation of the maximal eigenvalue $\lambda_{0}$ of $A$ and the minimal eigenvalue $\lambda_{\min}$. Let $\gamma:\Delta\rightarrow X$ be any holomorphic map from the disc to the generic fibre $X=\pi^{-1}(1)$. This defines a section
$\sigma_{\gamma}$ of $\cX$ by
$$  \sigma_{\gamma}(t)= (t, g(t) \gamma(t)) $$
with the usual interpretation of this formula at $t=0$. The map $g$ defines an isomorphism of the line bundle $\sigma_{\gamma}^{*}(\cL)$ with $\gamma^{*}(L)$ over the punctured disc. We define an integer $d(\gamma)$ to be the  first Chern class of $\sigma_{\gamma}^{*}(\cL)$ relative to this isomorphism. The next result immediately implies Lemma 3.
\begin{lem}
The minimal value of $d(\gamma)$, as $\gamma$ ranges over all such maps, is $-\lambda_{0}$, and the maximal value is $-\lambda_{\min}$.
\end{lem}

To prove this,   choose a lift $\tilde{\gamma}:\Delta\rightarrow V\setminus \{0\}$ of $\gamma:\Delta\rightarrow X$. Then $d(\gamma)$ is simply the order of the pole of the vector-valued function $g(t) \tilde{\gamma}(t)$ at $t=0$. It is then immediate that $-\lambda_{0}\leq d(\gamma)\leq -\lambda_{\min}$.
   Let $V_{-}\subset V$ be the subspace spanned by all eigenspaces
of $A$ with $\lambda_{\alpha}<\lambda_{0}$, and let $V^{+}$ be the subspace
spanned by all eigenspaces with $\lambda_{\alpha}>\lambda)_{\min}$. Choose any $\gamma_{+}$ with $\gamma_{+}(0)$ not in $\bP(V_{+})$. (This is possible because $X$ does not lie in $\bP(V^{+})$.) Then it is clear that 
$d(\gamma_{+})=-\lambda_{\min}$. In the other direction, we can suppose (for clarity) that $\cX$ is normalised so $\lambda_{0}=0$ and all $\lambda_{\alpha}\leq 0$. Choose  a smooth point $[y]$ of the central fibre $X_{0}$ which does not lie in $\bP(V_{-})$, that is the vector $y$ has a non-zero component in the zero eigenspace.  (This is possible since $X_{0}$ does not lie in $\bP(V_{-})$.)
  Choose a section $(t,\Gamma(t)$ of $\cX$ with $\Gamma(0)= y$. Then define
  $\tilde{\gamma_{-}}(t)= g(t)^{-1} \Gamma(t)= R(t)^{-1} t^{-A} \Gamma(t)$. Both
  $\tilde{\gamma_{-}}(t), g(t) \tilde{\gamma_{-}}(t)$ are holomorphic and non-vanishing at $t=0$ so  $d(\gamma_{-})=0$.

 {\bf Remark} More generally, one can give an intrinsic characterisation of  all the eigenvalues $\lambda_{\alpha}$ of $A$ via a filtration of the vector bundle $\pi_{*}(\cL)$ over $\Delta$, related to Szekelyhidi's filtration of the co-ordinate ring $\bigoplus H^{0}(X, L^{k})$. The matrix entries $R_{\alpha \beta}$ then appear, ore intrinsically, as representatives for the $H^{1}$ classes define the successive extensions of this filtration, relative to the trivialisation over $\Delta^{*}$.

 \subsection{Proof of Proposition 10 }
 We want to generalise the definitions of the Chow and Futaki invariants
to a \lq\lq relative'' case. Suppose we have families $\cX_{1}, \cX_{2}$ over
the disc and an isomorphism $\Psi$ of their restrictions to the punctured disc
$\Delta^{*}$. We suppose we have line bundles $\cL_{1},\cL_{2}$, very ample on
the fibres, and that $\Psi$ has a lift to the line bundles. Then the direct
images of $\cL_{1}^{k}, \cL_{2}^{k}$ are vector bundles over the disc with
an induced isomorphism between them over the punctured disc. Thus there is
a well-defined integer degree $-\tau_{k}(\cX_{1},\cX_{2})$ (where we suppress $\Psi$
in the notation).   If $\cX$ is a degeneration  as before, trivialised over $\Delta^{*}$, and if $\cX_{0}$
is the trivial family $X\times \Delta$ then 
$\tau_{k}(\cX,\cX_{0})= \tau_{k}(\cX)$.

 We assume that $\tau_{k}(\cX_{1}, \cX_{2})$ is, for large $k$, a polynomial of degree at most
$n+1$ in $k$. It seems clear that this is true as a general fact. In our particular
application we will be able to prove it by a separate argument, which we
give below, but to avoid complication we assume this for the time being.
Then we can define the Chow and Futaki invariants just as before. Thus
if $\tau_{k}(\cX_{1},\cX_{2})= I k^{n+1}+ O(k^{n})$ we let
$$  \Ch_{k}(\cX_{1},\cX_{2})= \frac{\tau_{k}}{N_{k}+1}- k\frac{I}{V}, $$
and $$\tFut(\cX_{1},\cX_{2})= \lim_{k\rightarrow \infty} \Ch_{k}(\cX_{1},\cX_{2}).$$

 Let $\cX$ be a non-trivial  degeneration, with general fibre $X=\pi^{-1}(1)$ as in the statement of Proposition 10. We can suppose that $\cX$ is normalsed. Write $\cX_{0}$ for the trivial degeneration $X\times \Delta$ and $\cL_{0}\rightarrow \cX_{0}$ for the pull-back of $L\rightarrow X$.
Let $\tilde{\gamma}_{-}:\Delta\rightarrow V\setminus\{0\}$ be a holomorphic map as considered in the proof of Lemma 8, so there is an induced map $\gamma_{-}$ from $\Delta$ to $X$. We regard this as  a section $\sigma_{2}:\Delta\rightarrow \cX_{0}$ with image a curve $\Sigma_{2}\subset \cX_{0}$. Let $\sigma_{1}$ be the section $\sigma_{\gamma_{-}}$ of $\cX$, as in the proof of Lemma 8, and write $\Sigma_{1}\subset \cX$ for the image. As in the proof of Lemma 8, the map $g$ induced an isomorphism of the line bundles
$\sigma_{1}^{*}(\cL), \sigma_{2}^{*}(\cL_{0})$ over the  disc $\Delta$.

  Let $\cX_{1}$ be the blow-up of $\cX$ along $\Sigma_{1}$ and $\cX_{2}$ be the blow-up of $\cX_{0}$ along $\Sigma_{2}$. Write $\cL_{E}$ for the line bundle defined by the exceptional divisor in both cases. For sufficiently large $\mu$  the line bundles $\cL^{\mu}\otimes \cL_{E}$ and $\cL_{0}^{\mu}\otimes \cL_{E}$ are ample on the fibres, and we write $\cX_{1,\mu}, \cX_{2,\mu}$ for the degenerations with these choices of polarising line bundles.
Thus we have invariants $\tFut(\cX_{1,\mu},\cX_{2,\mu}), \tFut(\cX,\cX_{0})$.
Our proof of Proposition 10 breaks into two parts. 
\begin{prop}
If $\mu$ is sufficiently large then 
$$\tFut(\cX_{1,\mu},\cX_{2,\mu})<\tFut(\cX,\cX_{0})= \tFut(\cX). $$
\end{prop}
\begin{prop}
If $X$ is bi-asymptotically stable then, for sufficiently large $\mu$,
$\tFut(\cX_{1,\mu}, \cX_{2,\mu})\geq 0$.
\end{prop}

We will now return to the issue arising in  the definition of the \lq\lq relative''
Futaki invariants. Suppose we can find a compact Riemann surface $B$ and
an embedding of the disc $\Delta$ in $B$ such that the family $\cX_{2}$ over the disc
extends to a family $\cY_{2}$ over $B$. Then, using the isomorphism between the families over the punctured disc,  we can also extend $\cX_{1}$ over
$B$ to get another family $\cY_{1}$, isomorphic to $\cY_{2}$ outside $\Delta$.
We define $\tau_{k}(\cY_{1}), \tau_{k}(\cY_{2})$ in the obvious way, using the degree
of direct images, and we clearly have
$$    \tau_{k}(\cX_{1}, \cX_{2})= \tau_{k}(\cY_{1})-\tau_{k}(\cY_{2}). $$
Then in such a situation we do know that $\tau_{k}(\cX_{1},\cX_{2})$ is a polynomial
of at most degree $n+1$, for large $k$, from the Grothendiek-Riemann-Roch
theorem applied to $\cY_{1}, \cY_{2}$. We can also define $\Ch_{k}(\cY_{i}), \tFut(\cY_{i})$
in an obvious way, so that 
$$   \tFut(\cX_{1}, \cX_{2})= \tFut(\cY_{1})-\tFut(\cY_{2}). $$

In our  situation we have a holomorphic map $\gamma_{-}:\Delta\rightarrow X$.
It is clear that any other map which agrees with $\gamma_{-}$ to a sufficiently high
order at $0$ will do just as well. Thus we need the following Lemma
\begin{lem}
Given any holomorphic map $\gamma_{-}:\Delta\rightarrow X$ and any integer
$s$ there is a holomorphic map $\gamma':\Delta\rightarrow X$ which
is equal to $\gamma_{-}$ to order $s$ at $0\in \Delta$ and which factors as
$\gamma'= j\circ \iota$ where $\iota:\Delta\rightarrow B$ is an embedding
in a compact Riemann surface $B$ and $j:B\rightarrow X$ is holomorphic.
\end{lem}
To see this choose a generic projection $p:X\rightarrow \bC\bP^{n}$ which
restricts to a holomorphic isomorphism on a neighbourhood of $\gamma(0)$.
We can approximate
$p\circ\gamma_{-}$ to arbitrarily high order at $0$ by a rational curve $j_{0}:\bC\bP^{1}\rightarrow
\bC\bP^{n}$. We can also suppose that $j_{0}$ is in general position relative
to the branch set of $p$. Then let $B$ be the normalisation of the curve
$p^{-1}(j_{0}(\bC\bP^{1}))$.

Given this Lemma, we can suppose that $\gamma_{-}$ itself extends to a Riemann surface $B$. We define $\cY_{0}$ to be the product family $X\times B$ and $\cY_{1}$ to be the family over $B$ equal to the product away from $\Delta$ and equal to $\cX$ over $\Delta$, as above.  We will now write $\Sigma_{2, B}$ for the graph of $j$ in $B\times X= \cY_{0}$ and $\Sigma_{1,B}\subset \cY_{1}$ for the section corresponding to $\Sigma_{1}$ over $\Delta$ and to the graph of $j$ away from $\Delta$.  We define $\cY_{1,\mu}, \cY_{2,\mu}$ by blowing up $\Sigma_{1,B}, \Sigma_{2,B}$  respectively.

We will now prove Proposition 12. What we need to show is that, if $X$ is
bi-asymptotically stable and  $\mu$ is sufficiently large, then
\begin{equation} \tFut(\cY_{1,\mu})\geq \tFut(\cY_{2,\mu}). \end{equation}
In fact we show that if $X$ is $(\mu,r)$ stable then $\Ch_{r}(\cY_{2,\mu})\geq
\Ch_{r}(\cY_{1,\mu})$.

    Recall the Chow construction, discussed in 3.1 above. Applied to a family $\cY\rightarrow B$ we get a line bundle $\lambda_{\Chow, \cY}$ over $B$ and we have
\begin{equation} c_{1}(\lambda_{\Chow, \cY})= \pi_{*}((c_{1} \cL)^{n+1}).
 \end{equation}
 Now consider our two families 
$\cY_{1,\mu}, \cY_{2,\mu}$ which are isomorphic away from the origin in $\Delta \subset
B$. Choose  trivialisations of $\pi_{*}(\cL_{\cY_{1,\mu}}), \pi_{*}(\cL_{\cY_{2,\mu}})$
The fibres of $\cY_{1,\mu},\cY_{2,\mu}$ over $\Delta$ can then be viewed as two families  of varieties in a fixed $N$-dimensional projective space
$ \bP(V)$, parametrised by $\Delta$, and we have two families of Chow points
in a projective space $\bP(W)$. Choose a holomorphic map $\sigma:\Delta\rightarrow
W\setminus\{0\}$ representing the Chow points of the fibres of $\cY_{2,\mu}$. Away from $0\in \Delta$ the
fibres of $\cY_{1,\mu}$ are isomorphic to those of $\cY_{2,\mu}$ and we get a meromorphic matrix valued function
$g(t)$ as in (3.1) and (5.2) above, such that the Chow point of the fibre of $\cY_{1,\mu}$ over $t$ is $[g(t) \sigma(t)]$. Let  $a$ be the order of the pole of $g(t) \sigma(t)$ and $b$  the order of the pole of $\det g$. This is a slight generalisation of the discussion in (3.1) because we are considering a family of orbits, varying with $t$.  But it is clear that if the orbits are stable then we have the same inequality 
$$\frac{a}{V}-\frac{b}{N+1}\geq 0 $$ On the other hand we can write
$$  a= c_{1}(\lambda_{\Chow, \cY_{1,\mu}})-c_{1}(\lambda_{\Chow, \cY_{2,\mu}})\ \ , b=  c_{1}(\pi_{*}
\cL_{Y_{1,\mu}}) - c_{1}(\pi_{*}\cL_{Y_{2,\mu}}). $$
So the inequality $a/V- b/(N+1)\geq 0$ is exactly the inequality $\Ch(\cY_{1,\mu})\geq \Ch(\cY_{2,\mu})$ that we want.

Now we turn to the proof of Proposition 11. We consider a general situation where we have
a family $\cL\rightarrow \cZ\rightarrow B$ and a section $\Sigma\subset \cZ$
contained in the smooth part of $\cZ$. We
blow up $\cZ$ along $\Sigma$ to obtain another family $\hat{\cZ}_{\mu}\rightarrow
B$ with an ample line bundle defined by $\mu$, as above. 
 Let $w$ be the degree of $\cL$ restricted
to $\Sigma\subset \cZ$. Then we have
\begin{prop}
In this situation
$$  \tFut(\hat{\cZ}_{\gamma})= \tFut(Z)- (w-\frac{I(\cZ)}{V})\mu^{1-n}+ O(\mu^{-n}) $$
\end{prop}

The proof is a straightforward calculation, similar to that of Stoppa but simpler because the blow-up locus lies in the smooth part.(One can also use the cohomological formula of Odaka \cite{kn:Odaka1}.)

We apply this formula to the pair $\cY_{0},\cY_{2,\mu}$  and $\cY, \cY_{1,\mu}$. Recall that we are assuming that our original family $\cX$ is normalised and the map $\gamma_{-}$ is chosen so that the relative first Chern class is zero. This implies that the degrees of the polarising line bundles over $\Sigma_{1,B}, \Sigma_{2,B}$ are equal. So Proposition 13 gives
$$  \left(\tFut(\cY_{1,\mu})-\tFut(\cY)\right)-
\left( \tFut(\cY_{2,\mu})-\tFut(\cY_{0})\right)= \frac{I(\cY)-I(\cY_{0})}{V} \mu^{1-n}+ O(\mu^{-n}). $$

The left hand side is $\tFut(\cX_{1,\mu},\cX_{2,\mu})- \tFut(\cX,\cX_{0})$. It is clear that $I(\cY)-I(\cY_{0})= I(\cX)$, so we have
$$\tFut(\cX_{1,\mu},\cX_{2,\mu})- \tFut(\cX,\cX_{0})= I(\cX)\mu^{1-n} + O(\mu^{-n}), $$ and Proposition 11 follows from Lemma 7.

\subsection{Proof of Proposition 1}

Recall that we have a Chow variety ${\rm Chow} \subset P(W)$ where $GL(N+1)$ acts linearly on $W$. For each point $s\in S$ we have an orbit $\Gamma_{s}$ in $\Chow$. Define  $d(s)$ to be the degree of the closure of $\Gamma_{s}$.
General facts of algebraic geometry imply that $d$ is a semi-continuous function so that the subset of $S$ on which takes its maximum value is Zariski open. Further, for a point $\sigma$ in this subset, the first condition of Definition 3  holds. Thus, replacing $S$ by this subset, we can suppose that this condition holds throughout $S$ (i.e. there is no \lq\lq splitting of orbits'').

Now we turn to the second condition. 
By restricting to a Zariski open subset of $S$ we can suppose that there is  a map $s\mapsto w_{s}\in W$ which chooses
a representative $[w_{s}]\in \Gamma_{s}\subset \Chow\subset \bP(W)$. 
Any point in the closure of an orbit $\Gamma_{s}$ can be written as the limit
as $t\rightarrow 0$ of $[L \Lambda(t)R(t)w_{s}]$ where $L$ is in $GL(N+1)$, $\Lambda(t)=t^{A}$ for a diagonal matrix $A$ and $R(t)$ is holomorphic and invertible across $t=0$. For each choice $c$ of diagonal entries of $A$ (i.e. the conjugacy class of the 1-parameter subgroup $\Lambda$), we can  suppose that $R$ is a lower triangular matrix whose entries are polynomials with degree bounded as in Lemma 6. Thus we can construct an algebraic variety $P_{c}$ and for each point $p\in P_{c}$ a map $g_{p}:\Delta^{*}\rightarrow GL(N+1, \bC)$ so that any point in the closure of $\Gamma_{s}$ arises as the limit when $t$ tends to $0$ of $[g_{p}(t)(w_{s})]$ for some $c$ and some
$p\in P_{c}$. We define a
 map of sets $E_{c}:S\times P_{c}\rightarrow \Chow$  by $E_{c}(s,p)= \lim_{t\rightarrow 0}[g_{p}(t)(w_{s})]$. We also define  an integer-valued function $\nu(s,p)$
   as the order of the pole of  the vector valued function $g_{p}(t)w_{s}$ at $t=0$.  Clearly $\nu$ is a bounded function. The limit $E_{c}(s,p)$ is not a continuous function of $(s,p)$ in general but it is so on any subset on which  $\nu$ is constant. Let $\Sigma_{\nu,c}\subset S\times P_{c}$ be the set where $\nu(s,p)=\nu$. The image of the map $E_{c}$ from $\Sigma_{\nu,c}$ to $\oJ$ is a quasi-projective variety and the union of these as $c$ runs over all the conjugacy classes is the whole of $\oJ$. Since there are only a countable number of conjugacy classes we can in fact choose a finite set of conjugacy classes $c$ which suffices to cover $\oJ$ in this way. Now consider the projection map $\pi_{S}$ from $\Sigma_{\nu,c}$ to $S$. By passing to a Zariski-open subset we may suppose that, for each $\nu, c$ in question, these are submersions mapping onto $S$. We claim then that in this situation the subset $GD$ is in fact  all of $S$. For given any pair $([Y],\sigma)$ in $\oJ$ we can find one of our finite set of pairs $(\nu,c)$ so that there is a $(\sigma,p)$ in $\Sigma_{\nu,c}$ with $[Y]=\lim_{t\rightarrow 0} [g_{p}(t) w_{\sigma}]$. Let $\pi_{S}^{-1}:N\rightarrow
\Sigma_{\nu,c}$ be a local right inverse to the projection map $\pi_{S}$ restricted to $\Sigma_{\nu,c}$ near $p$, with $\pi^{-1}(\sigma)=p$. For $s\in N$ write $g(s,t)= g_{\pi_{S}^{-1}(s)}(t)$. Then $g(s,t)$ defines a family of degenerations parametrised by $N$ as required.

\section{Proof of Proposition 7}

Let $\Lambda\rightarrow Z$ be a holomorphic Hermitian line bundle over a compact complex $n$-manifold with curvature $-i F_{\Lambda}$, where $F_{\Lambda}>0$.
We write $\rho_{\Lambda}$ for the function on $Z$ defined by the sum of squares of an orthonormal basis for the holomorphic sections using the volume form
$(n!)^{-1} F_{\Lambda}^{n}$. 

\begin{prop}
For each $n$ there   are constants $\kappa_{n}, K_{n}$ with the following effect. Let $(Z,\Omega)$ be any compact K\"ahler manifold of complex dimension $n$ and let $\Lambda\rightarrow Z$ be a holomorphic Hermitian line bundle with curvature $-i\kappa\Omega$. Let $B\subset Z$ be the metric ball of radius $1$ centred at $p\in Z$. Suppose that
\begin{enumerate}
\item The injectivity radius at $p$ is greater than $1$ (so the exponential map is invertible on $B$).
\item  The curvature bound $\vert \Riem\vert \leq 1 $ hold in $B$.
\item The metric has constant scalar curvature $S$ in $B$.
\item There is a Ricci curvature bound $\vert Ricci \vert \leq 1$ on all of $Z$.
\end{enumerate}
Then if $\kappa\geq \kappa_{0}$ we have
\begin{equation} \rho_{k}= (2\pi)^{-n}(1 + \kappa^{-1} S + E \kappa^{-2}), \end{equation}
where $$ \vert E\vert, \vert \nabla E\vert, \vert \nabla\nabla E \vert \leq K_{n}$$  in the ball of radius $1/2$ centred at $p$.
\end{prop} 

Proposition 7 can be deduced from Proposition 14  by rescaling. Let $\omega$ be a constant scalar curvature metric on $X$, as considered in 3.2, and let $x$ be a point in $\Omega_{r}\subset X$, with $r\leq 1$.  Rescale the metric by a factor $\beta r^{-2}$ where $\beta\geq 1$ will be chosen presently. The definition of $\Omega_{r}$ implies
that condition (2) in Proposition $14$ is satisfied by the rescaled metric (with $p=x$). We choose $\beta\geq C$ where $C$ is the constant in the definition of Property R, so condition (4) is satisfied. The diameter and Ricci bounds in the definition of Property R yield a uniform lower bound on the volumes of metric balls which in turn leads to a local injectivity radius bound in terms of the curvature. This means that we can choose $\beta$ so that the condition (1) holds. Now $\beta$ is fixed. We take $\Lambda=L^{k}$ so $\kappa = k r^{2}\beta^{-1}$. We take the constant $b_{0}$ in Proposition 7 to be $\beta \kappa_{0}$ which ensures that $\kappa\geq \kappa_{0}$. The function $\eta_{k}$ in Proposition 7 is $\kappa^{-2} E= k^{-2} r^{4} \beta^{2} E$ and the bounds asserted in Proposition 7 follow from those in Proposition 14, taking account of the rescaling of derivatives.

The proof of Proposition 14 is implicit in the literature, particularly Lu's paper \cite{kn:Lu}.  The key point is that in the standard proofs of the asymptotics the only control needed away from the ball $B$ is  the Ricci curvature bound. We will review a proof, but this is purely expository. See also the exposition in \cite{kn:BBS}.

  Generalising the $O(\ ), o(\ )$ notation
we write $\epsilon(\kappa)$ for any term  which is bounded by $C \kappa^{-m}$ for all
$m$. Essentially this means exponentially decaying terms, invisible in any
asymptotic expansion in inverse powers of $\kappa$.

To bring out the main point consider first a case when the metric is actually
{\it flat} in the embedded unit ball $B\subset Z$. We identify this
with the standard ball in $\bC^{n}$. We fix an identification of the fibre
of $\Lambda$ over the origin with $\bC$.  Rescale
the metric $\Omega$ by a factor $\kappa^{-1}$, so we work in a large ball  $B_{\sqrt{\kappa}}\subset \bC^{n}$, and
we are operating with the standard holomorphic Hermitian line bundle with
curvature $(i/2) \sum dz_{a} d\overline{z}_{a}$. Over this ball we
have a standard holomorphic section  $\sigma_{0}$ with  $\vert \sigma_{0}
\vert=\exp(-\vert z \vert^{2}/4)$. We cut of this section near the boundary
of the ball to get a compactly-supported section $\sigma_{1}$ which is approximately
 holomorphic in that
$$  \Vert \overline{\partial} \sigma_{1} \Vert = \epsilon(\kappa). $$
We can use a spherically symmetric cut-off function to do this, so $\sigma_{1}$
is $U(n)$-invariant in an obvious sense. We also have $\Vert \sigma_{1}\Vert
=(2\pi)^{n} (1+\epsilon(\kappa)$ and evidently $\sigma_{1}(0)=1$. 

Now transplant $\sigma_{1}$ to a section of $\Lambda$ over our manifold $Z$.
We work with the rescaled metric so that $Z$ has large volume $O(\kappa^{n})$.
We project $\sigma_{1}$ to the space of holomorphic sections using the standard
H\"ormander technique which we now recall. Since the Ricci curvature of the
original manifold is bounded that of the rescaled manifold is $O(\kappa^{-1})$.
The Weitzenbock formula on $\Lambda$-valued $(0,1)$ forms (for the rescaled
metric) takes the shape
$$   \Delta_{\db}= \frac{1}{2} \nabla^{*}\nabla + {\rm Ric} + 1 $$
so once $\kappa>2$, say, we have $\Delta_{\db}\geq \frac{1}{2}$ say, and the inverse
operator $G$ has $L^{2}$-operator norm at most $2$. Now 
$$ \sigma_{2}=\sigma_{1}- \db^{*} G \db \sigma_{1} $$
is a holomorphic section. From the identity
$$  \Vert \db^{*} G \db \sigma_{1} \Vert^{2}= \langle G \db \sigma_{1}, \sigma_{1}
\rangle, $$
we get  
$$ \Vert \sigma_{2} - \sigma_{1} \Vert^{2} \leq 2 \Vert \sigma_{1} \Vert
\Vert \db \sigma_{1} \Vert = \epsilon(\kappa). $$

Go back to the unit ball (say) in $\bC^{n}$. Here the difference $\sigma_{2}-\sigma_{1}$
is holomorphic so the $L^{2}$ bound above gives a pointwise bound and
$$  \vert \sigma_{2}(0)-1 \vert = \epsilon(\kappa). $$

Let $\tau$ be any section over $Z$ which vanishes at the origin. The
inner product $\langle \sigma_{1}, \tau \rangle$ is an integral over the
ball $B_{\sqrt{\kappa}}$ and this obviously vanishes by symmetry (considering
the Taylor series of $\tau$ in our given trivialisation). Finally, for convenience,
set $\sigma=\sigma_{2}/\Vert \sigma_{2} \Vert$. What we have achieved is
a section $\sigma$ with the three properties
\begin{enumerate}
\item $ \Vert \sigma \Vert=1 ;$
\item $\vert \sigma(0)\vert^{2} = (2\pi)^{-n} (1 + \epsilon(\kappa));$
\item $\langle \sigma,\tau\rangle= \epsilon(\kappa) \Vert \tau \Vert$ if $\tau(0)=0$.
\end{enumerate}

No more  analytical input is required.  Let $\eta$ be the section representing
evaluation at $0$, i.e.
$$\langle \tau, \eta \rangle = \tau(0) $$
for all $\tau$. By definition the Bergman function at $0$ is $\rho_{\Lambda}(0)=
\Vert \eta \Vert^{2}$. What we need is
\begin{lem}
Let $\eta, \sigma$ be two elements of a Hilbert space such that
\begin{enumerate}
\item $\Vert \sigma \Vert=1$;
\item $\langle \sigma, \eta \rangle= 1+\epsilon(\kappa)$;
\item $\langle \sigma, \tau \rangle= \epsilon(\kappa) \Vert \tau \Vert $ for any
$\tau$ with $\langle \tau,\eta\rangle=0$.
\end{enumerate}
Then $\Vert \eta \Vert^{2} = 1+\epsilon(\kappa)$.
\end{lem}
The proof is an elementary exercise (which takes place in the plane spanned
by $\sigma, \eta$).In our case the three hypotheses are re-statements of
the properties above and we conclude that, in this locally flat situation,
 $\rho_{\Lambda}(0)= 1+ \epsilon(\kappa)$. The bounds on all derivatives of $\rho_{\Lambda}$ follow by a straightforward extension of the argument.

Now we go on to the general case. Note first that the constant scalar curvature condition gives estimates on all derivatives of the metric in the interior of $B$ (as in Section 4)---this is the only role of the condition in the statement of Proposition 14.  By the same argument as for the Lemma above,
if we produce a holomorphic section $\sigma$ with $\Vert \sigma \Vert=1$, with
$$  \sigma(0)= (2\pi)^{-n/2}(1+ A\kappa^{-1}+ O(\kappa^{-2}) , $$
and $$  \vert \langle \sigma, \tau\rangle\vert\leq C \kappa^{-1} \Vert \tau \Vert,
 $$
for all holomorphic sections $\tau$ vanishing at the origin, then $\rho_{\Lambda}(0)=
(2\pi)^{n}(1+2A\kappa^{-1})^{2}+ O(\kappa^{-2})$. Consider  holomorphic co-ordinates $w_{a}$ centred
on the given point in the manifold. A K\"ahler potential $\phi$ has a Taylor
series which we can obviously suppose begins as
$$   \phi(w)= \sum_{a} w_{a} \overline{w}_{a} + O(w^{3}) $$
Write the cubic term schematically as $(3,0)+(2,1)+(1,2)+(0,3)$ in terms
of the degree in $w_{a}, \overline{w}_{a}$. By a change of co-ordinates of
the form
$$\tilde{w}_{a}= w_{a}+ C_{a bc} w_{b} w_{c}$$
we can reduce to the case when the $(2,1)$ and $(1,2)$ terms (which are complex
conjugate) vanish. Then by adding the real part of a holomorphic quartic
function to the Kahler potential we can remove the $(3,0)$ and $(0,3)$ terms.
Similarly we can remove all the quartic terms in the Taylor expansion except
for $(2,2)$ and all the quintic terms except for $(2,3)+(3,2)$. So we can
suppose that the K\"ahler potential $\phi$ is
$$  \ \sum w_{a} \overline{w}_{a} + \sum P_{a b c d} w_{a} w_{b} \overline{w}_{c}\overline{w}_{d}
+ (\sum Q_{a b c d e}w_{a} w_{b} w_{c}\overline{w}_{d} \overline{w}_{e} +
{\rm cx.\  conjugate}) +O(w^{6}). $$

Next we rescale co-ordinates writing $w_{a}= \kappa^{-1/2} z_{a}$ and setting
$\Phi(z)= \kappa \phi(w)$. Thus 
$$  \Phi(z)= \vert z \vert^{2} + \kappa^{-1} P(z) + \kappa^{-3/2} Q(z) + O(\kappa^{-2}),$$
in an obvious notation. We work over a ball $\vert z\vert \leq R$ where
we can take $R$ to be a very small multiple of $\kappa^{1/4}$ so that $\Phi(z)-
\vert z \vert^{2}, \kappa^{-1} P(z), \kappa^{-3/2} Q(z) $ are all very small over the
ball. The volume form  $(\dbd \Phi)^{n}$ in these co-ordinate can be written
$$  J= 1+ \kappa^{-1} p(z) + \kappa^{-3/2} q(z) +O(\kappa^{-2}), $$
where for example
$$p(z)= 4 \sum_{a,b,c} P_{a b c b} z_{a}\overline{z}_{c}. $$

The choice of a Kahler potential precisely corresponds to the choice of a
local trivialisation of our line bundle and hence a local holomorphic section
$\sigma_{0}$ with $\vert \sigma_{0}\vert^{2}= e^{-\Phi}$. Just as before
 we can modify $\sigma_{0}$ to get a global holomorphic section and this
only introduced terms which are $\epsilon(\kappa)$ so which we can ignore. Regard $\sigma_{0}$
as a discontinuous section of the line bundle over the whole manifold, extending
by zero outside our ball. We want to show first that 
$$  \vert \langle \sigma_{0}, \tau\rangle \vert \leq C k^{-1}, $$
for all holomorphic sections vanishing at the origin. Second, we want to find
a number $a$ such that $  \Vert \sigma_{0} \Vert^{2} = 1+ a\kappa^{-1} + O(\kappa^{-2})$.
Then we will have established what we need with $A=-a/2$ (since $\sigma_{0}(0)=1$
by construction). 

Now $$  \vert \langle \sigma_{0}, \tau\rangle \vert \leq C \kappa^{-1} \int_{\vert
z\vert \leq R} (1+ \vert z\vert^{4}) e^{-\vert z\vert^{2}} \tau(z), $$
where $\tau(z)$ is the representative of $\tau$ in our local trivialisation.
Here we use the fact that
$$  \int_{\vert z\vert \leq R} \tau(z) e^{-\vert z\vert^{2}}= 0$$
when $\tau(0)=0$. We obtain the desired estimate using the Cauchy-Schwartz
inequality in the weighted norm. So it just remains to compute $\Vert \sigma_{0}
\Vert^{2}$ which is
$$   \int_{\vert z\vert \leq R} e^{-\vert z\vert^{2}} (1+ \kappa^{-1} P(z)+ \kappa^{-3/2}
Q(z)+O(\kappa^{-2})(1+ \kappa^{-1} p(z) + \kappa^{-3/2} q(z) +O(\kappa^{-2}). $$
Here we have used the Taylor series to expand the exponential term $e^{-\Phi}$
and we have skipped over some rather routine estimates. 

If $z^{I}, \overline{z}^{J}$ are any monomials such that
$$  \int_{\vert z\vert \leq R} z^{I} \overline{z}^{J} e^{-\vert z\vert^{2}}
$$
is non zero then we must have $\vert I\vert=\vert J\vert$. (To see this,
consider the action of multiplication by $e^{i\theta}$ .) It follows that
the integrals appearing in the $\kappa^{-3/2}$ terms above vanish. Extending the
range of integration introduces errors $\epsilon(\kappa)$ so we get 
$$ \Vert \sigma_{0}\Vert^{2}= \int_{\bC^{n}} \left(1+ \kappa^{-1}(P(z)+p(z))\right)e^{-\vert
z\vert^{2}} +O(\kappa^{-2}). $$
The integral here is straightforward  to calculate. We can also argue as
follows. The tensor $P_{a b c d}$ is, from an invariant point of view, an
element of $s^{2}(V) \otimes s^{2}(V)^{*}$ where $V$ is the cotangent space.
The symmetric power $s^{2}(V)$ is an irreducible representation of $U(n)$
so there is, up to a multiple, just one $U(n)$-invariant contraction $s^{2}(V)
\otimes s^{2}(V)^{*}\rightarrow \bC$. This is given by $$c(P)= \sum_{a,b}
P_{a b a b}. $$
It is clear that the scalar curvature and the integrals appearing in the
$O(\kappa^{-1})$ term above are both invariant contractions of $P$ hence multiples
of $c(P)$. This argument shows that the co-efficient we are after must be
some universal multiple of the scalar curvature and of course we can identify
the multiple from the Hirzebruch-Riemann-Roch formula.




\begin{thebibliography}{99}
\bibitem{kn:AP} Arrezzo, C. and Pacard, F. {\em Blowing up and desingularizing constant scalar curvature K\"ahler metrics} Acta Math. 196 179-228 2006
\bibitem{kn:BBS} Berman, R., Berndtsson, B., and Sj\"ostrand, J. {\em A direct approach to Bergman kernal asymptotics for positive line bundles} Ark. Math. 46 197-217 2008
\bibitem{kn:Caf} Cafarelli, L. {\em Interior a priori estimates for solutions of fully nonlinear equations} Annals of Math. 130 189-213 1989
\bibitem{kn:CheegerNaber} Cheeger, J. and Naber, A. {\em Lower bounds on Ricci curvature and quantative behaviour of singular sets} Inventiones Math. 321-339 2013
\bibitem{kn:CD} Chen, X. and Donaldson, S. {\em Volume estimates for K\"ahler-Einstein metrics and rigidity of complex structures} Jour. Differential Geometry  93 191-201 2013
\bibitem{kn:CDS} Chen, X., Donaldson, S. and Sun, S. {\em K\"ahler-Einstein metrics and stability} Int. Math. Res. Notices No 8 2110-2125 2014
\bibitem{kn:CH} Chen, X. and He, W. {\em On the Calabi flow} American Jour. Math. 130 539-570 2008
\bibitem{kn:CL} Cheng, S. and Li, P. {\em Heat kernel estimates and lower bounds of eigenvalues} Comm. Math. Helveticii 56 379-91 1976 
\bibitem{kn:D1} Donaldson, S. {\em Scalar curvature and projective embeddings, I} Jour Differential Geom. 59 479-522 2001 
\bibitem{kn:D2} Donaldson, S. {\em  K\"ahler geometry on toric manifolds, and some other manifolds with large symmetry} In: Handbook of Geometric Analysis, No. I Advanced Lectures in Mathematics Vol 7 International Press 2008 29-75
\bibitem{kn:D3} Donaldson,S. {\em Discussion of   the K\"ahler-Einstein problem} http://wwwf.imperial.ac.uk/$\sim$skdona/KENOTES.PDF
\bibitem{kn:D4} Donaldson, S. {\em Stability, birational transformations and the K\"ahler-Einstein problem}  Surveys in Differential Geometry Vol XVII
 Int. Press Boston MA 2012
\bibitem{kn:DS}  Donaldson, S. and Sun, S. {\em Gromov-Hausdorff limits of K\"ahler metrics and algebraic geometry} Acta Math. 213 63-106 2014
\bibitem{kn:JMR} Jeffres, T., Mazzeo, R. and Rubinstein, Y. {\em K\"ahler-Einstein metrics with edge singularities} arxiv 1105.5216
\bibitem{kn:LX} Li, C. and Xu, C. {\em Special test configurations and K-stability of Fano manifolds} Annals of Math. 180 197-232 2014
\bibitem{kn:KM} Knudsen, F. and Mumford, D. {\em Projectivity of moduli spaces of stable curves, I. Preliminaries on \lq\lq det'' and \lq\lq div''} Math. Scand. 39 19-55 1976
\bibitem{kn:Lu}  Lu, Z. {\em On the lower order terms in the asymptotic expansion of Tian-Yau-Zelditch} American Jour. Math. 122 235-273 2000
\bibitem{kn:NW} Nijenhuis, A. and Wolf. W. {\em Some integration problems
in almost-complex and complex manifolds} Annals of Math. 77 424-489 1963
\bibitem{kn:Odaka1}  Odaka, Y. {\em The Calabi conjecture and K-stability} Int. Math. Res. Notices 10 2272-2288 2012
\bibitem{kn:Odaka2} Odaka, Y. {\em On the moduli of K\"ahler-Einstein Fano manifolds} 1211.4833
\bibitem{kn:PS} Phong, D. and Sturm, J. {\em Scalar curvature, moment maps and the Deligne pairing} American J. Math. 126 693-712 2004
\bibitem{kn:Stoppa} Stoppa, J. {\em K-stability of constant scalar curvature K\"ahler manifolds } Advances in Math. 221 1397-1408 2009
\bibitem{kn:T} Tian, G. {\em K\"ahler-Einstein metrics with positive scalar curvature} Inventiones Math. 13 1-37 1997
\bibitem{kn:Gabor} Szekelyhidi, G. {\em Filtrations and test configurations} arxiv 1111.4986
\bibitem{kn:Zhang} Zhang, S. {\em Heights and reductions of semistable varieties}
Compositio Math. 104 77-105 1996 \end{thebibliography}
\end{document}